\documentclass[a4paper, 10pt]{amsart} 

\theoremstyle{plain} \newtheorem{thm}{Theorem}[section]
   \newtheorem{lem}[thm]{Lemma}
  \newtheorem{cor}[thm]{Corollary} 
  \theoremstyle{definition} 
   
  \theoremstyle{remark} \newtheorem{rem}[thm]{Remark}
  \numberwithin{equation}{section}

\def\oc{\mathcal{O}} \def\oce{\mathcal{O}_E} 
\def\ac{\mathcal{A}} \def\rc{\mathcal{R}} \def\mc{\mathcal{M}}

 \def\tly{{Y}} \def\om{\omega}

\def\cb{\mathbb{C}}  
\def\pb{\mathbb{P}} \def\pbe{\mathbb{P}(E)} 
\def\zbb{\mathbb{Z}}

\def\hb{\bold{H}} \def\fb{\bold{F}} \def\eb{\bold{E}}
\def\pbb{\bold{P}}

\def\nab{\overline{\nabla}} \def\n{|\!|} \def\spec{\textrm{Spec}\,}
\def\cinf{\mathcal{C}_\infty} \def\d{\partial}
\def\db{\overline{\partial}}
\def\hess{\sqrt{-1}\partial\!\overline{\partial}}
 \def\lra{\longrightarrow}
\begin{document}

\title{Hodge metrics and positivity of direct images}
\author{Christophe Mourougane, Shigeharu Takayama} \maketitle

{\textbf{Abstract.} }  \textit{Building on Fujita-Griffiths method of
computing metrics on Hodge bundles, we show that the direct image of
an adjoint semi-ample line bundle by a projective submersion has a
continuous metric with Griffiths semi-positive curvature.  This shows
that for every holomorphic semi-ample vector bundle $E$ on a 
complex manifold, and every positive integer $k$, the vector bundle
$S^kE\otimes\det E$ has a continuous metric with Griffiths
semi-positive curvature. If $E$ is ample on a projective manifold, 
the metric can be made smooth and Griffiths positive.}

%%%%%%%%%%%%%%%%%%%%%%%%%%%%%%%%%%%%%%%%%%%%%

\section{Introduction}
We study the positivity properties of direct images of adjoint line
bundles. We get

\begin{thm} \label{semi-ample}
Let $\phi~:~Y\lra X$ be a projective submersion between two complex
manifolds.  Let $L$ be a semi-ample line bundle on $Y$.  Then the
vector bundle $\phi_\star(K_{Y/X}\otimes L)$ has a continuous metric
with Griffiths semi-positive curvature.
\end{thm}

\begin{thm} \label{ample}
Let $\phi~:~Y\lra X$ be a submersion between two complex projective
manifolds.  Let $L$ be an ample line bundle on $Y$.  Then the vector
bundle $\phi_\star(K_{Y/X}\otimes L)$ has a smooth metric with
Griffiths positive curvature.
\end{thm}

 On grassmanian manifolds, abelian varieties and toric varieties,
these results were already proved either by the Castelnuovo-Mumford
criterion for global generation or by trying to mimic the Frobenius
morphisms over $\cb$ (\cite{ma}, \cite[Th{\'e}or{\`e}me 3]{mo}).  These
theorems should be compared with by now classical results due to
Fujita, Kawamata and others in algebraic geometry setting (see for
example~\cite[Chapter 6.3.E]{laz}).  These constructions may help to
find some topological properties of algebraic sub-varieties under
ampleness assumptions for the normal bundle (see for
example~\cite{ful} and~\cite{laz} part two).  The interesting outcome
is that metrics could be used to construct Morse functions, whose
indices are computed by the curvature.

For the main application, we consider a holomorphic vector bundle $E$
on a compact complex manifold, and intend to construct metrics on
vector bundles associated to $E$ which would reflect algebraic
positivity properties of~$E$.  Just note that on the variety
$\pi~:~\pbe\lra X$ of rank one quotients of $E$ the direct image
$\pi_\star(K_{\pbe/X}\otimes \oce (k+r))$, where $r$ is the rank of
$E$, is $S^kE\otimes\det E$ to infer
 
\begin{cor} 
Let $E$ be a semi-ample vector bundle on a complex manifold.
Then for all positive integer $k$, the vector bundles $S^kE \otimes
\det E$ have continuous metrics with Griffiths semi-positive
curvature.
\end{cor}

\begin{cor} 
Let $E$ be an ample vector bundle on a complex projective manifold.
Then for all positive integer $k$, the vector bundles $S^kE \otimes
\det E$ are Griffiths positive.
\end{cor}
The theory of resolution of the $\db$-equation with $L^2$-estimates
for example shows that Griffiths positivity and ampleness are
equivalent for line bundles. This implies through the curvature
computation of $\oce(1)$ that Griffiths positive vector bundles are
ample.  The converse is a problem raised by Griffiths~\cite[problem
(0.9)]{Gr}, and solved positively on curves by Umemura~\cite{um} using
the concept of stability (see also~\cite{campana}).  Our results
provide a weak answer to Griffiths's problem.  The appearance of the
determinant line bundle has the same origin than its appearance in the
vanishing theorem of Griffiths for the cohomology of ample vector
bundles.

The first idea in the proof is to use natural cyclic coverings to be
able to use Hodge metrics on the direct image of the relative
canonical sheaf under a proper K{\"a}hler surjective map between two complex
manifolds.  This may be seen as a metric aspect of Ramanujam's idea to
reduce vanishing theorem to topological properties~\cite{ra} (see
also~\cite{ko}).

Applying Griffiths method of computing metrics on Hodge bundles, we
will be able to compute the curvature of direct image of the structure
sheaf for proper K{\"a}hler submersions. It turns out that the metric on the
top direct image of the structure sheaf of the source manifold is
semi-negatively curved.

We then have to determine the singularities of the Hodge metrics.  For
 a proper surjective map $f~:~Z^{n+r}\lra X^n$, the Jacobian ideal
 $Jac_f$ is the ideal of $\oc_Z$ generated by the $n\times n$- minors
 of the matrix of the differential of $f$ (computed in any chart). The
 discriminant locus $\Delta_f \subset X$ is defined to be the image
 under $f$ of the sub-scheme of $Z$ where the map $f$ is not
 submersive, that is the sub-scheme cut out by $Jac_f$. Denote
 $X-\Delta_f$ by $X^0$ and $f^{-1}(X^0)$ by $Z^0$. The map $f^0
 :=f_{|Z^0}~:~Z^0\lra X^0$ is then a projective submersion.  With
 these notations, our main technical lemma generalizes Fujita's
 observation~\cite[lemma 1.12]{fu}.

\begin{lem}\label{mainlemma}
Let $f~:~Z\lra X$ be a proper K{\"a}hler surjective map between two
complex manifolds.
\begin{enumerate}
\item Then, the Hodge metric on $(f^0)_\star (K_{Z^0/X^0})$ extends as
  a metric with poles on $f_\star (K_{Z/X})$.
\item The Hodge metric on $(f^0)_\star (K_{Z^0/X^0})$ extends as a
  smooth metric on $f_\star (Jac_f \otimes K_{Z/X})$.
\end{enumerate}
\end{lem}
 This is far simpler than general results obtained by Kawamata,
Z{\"u}cker, Koll{\'a}r and Cattani-Kaplan-Schimd.

This paper is a revised and expanded version of our preprint
\cite{mt}.  Shortly before we ended the writing of our text, Bo
Berndtsson gave a nice proof of the Nakano positivity of the direct
image of adjoint ample line bundle in the same setting than
ours~\cite{be}.  We nevertheless feel worth to display our techniques
for these are different from his.

{\it Acknowledgment}. The first named author warmly thanks Indranil
Biswas for a collaboration in an attempt to prove similar results.

%%%%%%%%%%%%%%%%%%%%%%%%%%%%%%%%%%%%%%%%%%%%%

\section{The tools}

\subsection{Ampleness and positivity}
\label{ap}

We refer to~\cite{H-66} or~\cite[\S 6]{laz} for basics about ample
vector bundles and to~\cite{Gr} or~\cite[chapter VII]{de-book} for
basics about positive vector bundles.  All vector bundles are assumed
to be holomorphic.

A vector bundle $E$ on a compact complex manifold $X$ is said to be
{\it semi-ample}\/ if for some positive integer $k$, its symmetric
power $S^k E$ is generated by its global sections.  Associated to $E$,
we have $\pi~:~\pbe\lra~X$ the variety of rank one quotients of $E$
together with its tautological quotient line bundle $\oce (1)$.  The
semi-ampleness of $E$ is rephrased that for every $x \in X$, every
section $s_x$ of $\oce (k)$ over the fiber $\pb (E_x)$ extends to a
global section of $\oce (k)$ over $\pb (E)$.  This in particular
implies that $\oce (k)$ is generated by its global sections.  A vector
bundle $E$ is said to be {\it ample}\/ if its associated line bundle
$\oce (1)$ is ample on $\pbe$.  This in particular implies the
existence of an integer $k$ such that for every $x \in X$, every
section $s_x$ of $\oce(k)$ over the first infinitesimal neighborhood
of the fiber $\pb (E_x)$ extends to a global section of $\oce (k)$
over $\pb (E)$.

Recall the formula for the curvature of the Chern connection of the
quotient metric $h_q$ on $\oce (1)$ of a chosen metric $h$ on
$E$. Here $a^\star$ parametrizes rank one quotients of $E$ and
$FS(h_x)$ is the Fubini-Study metric of the hermitian form $h_x$ on
$E_x$.
\begin{eqnarray*}
%\label{courbure relative}
\nonumber \Theta (\oce (1),h_q)&=& \hess \log \n
e^\star_0+z_1e^\star_1+\cdots z_{r-1}e^\star_{r-1}\n^2\\ &=& \Theta
(\oc_{\pb (E_x)} (1), FS (h_x)) -\frac{\langle
\pi^\star\Theta(E^\star, h^\star) a^\star, a^\star\rangle} {\langle
a^\star,a^\star\rangle}.
\end{eqnarray*}

A vector bundle $E$ is said to be {\it Griffiths positive}, if it can
be endowed with a smooth hermitian metric $h$ such that for all $x\in
X$ and all non-zero decomposable tensors $v \otimes e\in TX_x\otimes
E_x$, the curvature term
$\langle\Theta(E,h)(v,\overline{v})e,e\rangle_h$ is positive, where
$\Theta(E,h)\in \mathcal{C}_{1,1}^{\infty}(X,Herm (E))$ is the
curvature of the Chern connection $\nabla_{E,h}$ of $(E,h)$.

A continuous hermitian metric $h$ on a vector bundle $b~:~E\lra X$ is
said to be {\it Griffiths positive}, if there exists a smooth positive
real $(1,1)$-form $\omega_X$ on $X$ such that in the sense of currents
$$ -\sqrt{-1}\d\!\db\log h(\xi) \geq b^\star \omega_X,
$$ where $h$ is seen as a continuous quadratic function on the total
space $E-X\times \{0\}$.  At the points where the metric $h$ is
smooth, these two notions of Griffiths positivity coincide.

\subsection{Cyclic covers}
\label{sec2}

A reference for this part is~\cite[\S~3]{es}.  Take a semi-ample line
bundle $L$ on a complex manifold $Y$ and fix a positive integer $k$
such that the tensor power $L^k$ is generated by its global sections.
Then Bertini's theorem~(see for example~\cite[ page 137]{gh}) insures
that a generic section $s$ of $L^k$ over $Y$ is transverse to the zero
section (i.e. $ds_{|D_s}~:TY_{|D_s}\lra L^k_{|D_s}$ is surjective),
and defines a smooth divisor $D_s:=(s=0)$.  Let
$$ p = p_s~: \tly_s\lra Y
$$ be the cyclic covering of $Y$ obtained by taking the $k$-th root
out of $D_s$.  The space
$$ \tly_s := \{l\in L / l^k=s(p(l)) \}
$$ is a smooth hypersurface of the total space $L$, and the map $p$ is
a finite cover of degree $k$ totally ramified along the zero locus
$D_s$ of $s$.  The space $\tly_s$ may also be described as the
spectrum $\spec \ac_s$ of the algebra
$$ \ac_s~:=\frac{\oplus_{i=0}^{+\infty} L^{-i}} {\left(l^\star
-\check{s} (l^\star)\ , \ \l^\star\in L^{-k}\right) }
$$ where $\check{s}$ is the sheaf inclusion $L^{-k}\stackrel{\times
s}\rightarrow\oc_{Y}$.  The direct image of the structure sheaf
$\oc_{\tly_s}$ is hence $p_\star \oc
_{\tly_s}=\ac_s\simeq\oplus_{i=0}^{k-1} L^{-i}$.

In our setting, we consider a projective submersion $\phi~:~Y\lra X$
of relative dimension $r$
 between two complex manifolds and $L$ a semi-ample line bundle on
 $Y$. A generic section of $L^k$ gives rise to a covering of $Y$ which
 we in fact regard as a family of coverings of the fibers of $\phi$.
 Because $p_s$ is a finite morphism, the spectral sequence of
 composition of direct image functors reduces to the following:
\begin{eqnarray*} 
%\label{iso}
\nonumber \rc ^{r} (\phi\circ p_s)_\star \oc_{\tly_s} &=&\rc ^{r}
\phi_\star \left(\rc ^0 p_\star \oc_{\tly_s} \right) \\ \nonumber
&=&\rc ^{r} \phi_\star \left(\oplus_{i=0}^{k-1} L^{-i}\right)
=\oplus_{i=0}^{k-1} \phi_\star\left(K_{Y/X}\otimes L^{i}\right)^\star.
\end{eqnarray*}
Here we have used Serre duality on the fibers of the smooth morphism
$\phi$.

\subsection{The Hodge metric}
We recall the basics on geometric variations of Hodge structures and
Griffiths's computations (\cite[theorem 6.2]{gr}) of the curvature the
Hodge metric (see also~\cite[$\S$ 7]{sc} and~\cite[chapter~10]{vo}).
Here we assume that $f : Y \lra B$ is a proper K{\"a}hler submersion of
relative dimension $r$, 
in particular we regard
$f : Y \lra B$ as a smooth family of compact K{\"a}hler 
manifolds of dimension $r$.  Fix a non-negative integer $d$.  The
local system $\rc^df_\star\cb$ can be realized as the sheaf of germs
of the flat sections of the holomorphic vector bundle $\hb^d_\cb$
associated with the locally free sheaf $(\rc^df_\star\cb)\otimes\oc_B$
endowed with the flat holomorphic connection~$\nabla~:\hb^d_\cb\lra
\Omega^1_B\otimes\hb^d_\cb$, the Gauss-Manin connection.  By
semi-continuity and Hodge decomposition, the vector spaces
$H^{p,d-p}(Y_b,\cb)$ ($b\in B$) have constant dimension. By elliptic
theory they hence form a differentiable sub-bundle $\hb^{p,d-p}$ of
$\hb^d_\cb$.  Denote by $\fb^p$ the differentiable sub-bundle $\oplus
_{i\geq p} \hb^{i,d-i}$ of $\hb^d_\cb$.  By a theorem of Griffiths,
the $\fb^p$ have natural structure of holomorphic sub-vector bundles
of $\hb^d_\cb$.  A relative Dolbeault theorem identifies $\eb^p:=
\fb^p/\fb^{p+1}$ with the holomorphic vector bundle associated with
the locally free sheaf $\rc^{d-p}f_\star\Omega_{Y/B}^p$.

We now recall the construction of the Hodge metric on the primitive
part of $\eb^p$. We fix a family $\eta_b$ ($b\in B$) of polarizations
given by a section of $\rc^2f_\star \zbb$ (e.g.\ the family of Chern
classes $c_1(L_{|Y_b})$ for an ample line bundle $L$ on $Y$).  By the
hard Lefschetz theorem, the bilinear form on the fibers of $\hb^d_\cb$
given by
$$S(c_1,c_2):=(-1)^{\frac{d(d-1)}{2}}\int_{Y_b}\eta_b^{r-d}\wedge
c_1\wedge c_2$$ is non-degenerate. We define the primitive cohomology
to be $\pbb^d:=Ker (\eta^{r-d+1}~: \hb^d \lra \hb^{2r-d+2})$, which is
also a differentiable sub-bundle of $\hb^d$.  By the Hodge-Riemann
bilinear relations, the differentiable sub-bundles $\hb^{p,d-p}$ and
$\hb^{p',d-p'}$ are orthogonal unless $p+p'=d$ and
$$h(c):=(\sqrt{-1})^{p-q} S(c,\overline{c})$$ defines a positive
definite metric on $\hb_{prim}^{p,d-p}:=\hb^{p,d-p}\cap \pbb^d$. We
set $\fb^p_{prim}:=\fb^p\cap \pbb^d$. Those bundle also have natural
holomorphic structures.  We also set
$\eb^p_{prim}:=\fb^p_{prim}/\fb^{p+1}_{prim}$.  The fiber-wise
isomorphism of $(\eb^p_{prim})_b$ with
$H_{prim}^{p,d-p}(Y_b,\cb)\subset H^d(Y_b,\cb)$ enables to equip the
holomorphic vector bundle $\eb_{prim}^p$ with a smooth positive
definite hermitian metric, called the Hodge metric.

We need some definitions in order to express the curvature of the
corresponding Chern connection.  First recall the transversality
property of the Gauss-Manin connection $\nabla
\fb^p\subset\Omega^1_B\otimes \fb^{p-1}$ which accounts for the
Cartan-Lie formula of the derivative of a family of cohomology classes
(see~\cite[proposition 9.14]{vo}).  Denote by
$\nab^p~:\eb^p\lra\Omega^1_B\otimes\eb^{p-1}$ the $\oc_B$-linear map
built by first lifting to $\fb^p$ applying the Gauss-Manin connection
and projecting to $\eb^{p-1}$.  The second fundamental form in
$\cinf^{1,0}(B,Hom(\fb^p, \hb^d_\cb/\fb^p))$ of the sequence
$$ 0 \lra \fb^p \lra \hb^d_\cb \lra \hb^d_\cb/\fb^p \lra 0
$$ with respect to the Gauss-Manin connection (or equivalently with
the flat metric on $\hb^d_\cb$) actually induces
$\nab^p~:\eb^p\lra\Omega^1_B\otimes\eb^{p-1}$.  Formulas for the
curvature of quotient hermitian holomorphic vector bundles then lead
to

\begin{thm} \label{Gr1} \cite[theorem 5.2]{gr}\,
The curvature $\Theta (\eb_{prim}^p)$ of the holomorphic vector
bundle~$\eb_{prim}^p$ endowed with its Hodge metric is given by
\begin{eqnarray*}
\langle \Theta (\eb_{prim}^p)(V,\overline{V})\sigma,\sigma
\rangle_{Hodge} = \langle \nab^p_V \sigma, \nab^p_V \sigma
\rangle_{Hodge} -\langle(\nab^{p+1}_V)^\star
\sigma,(\nab^{p+1}_V)^\star \sigma \rangle_{Hodge}
\end{eqnarray*}
where $V$ is a local vector field on $B$ and $\sigma$ a local section
of $\eb^p$.
\end{thm}

We now apply this result in the case of the family of the cyclic
covers $\phi\circ p_s~: {\tly_s} \lra X$ obtained by taking the $k$-th
root of a section $s$ of $L^k$ transverse to the zero section as in \S
~\ref{sec2}.  We have to restrict the study over Zariski open sets
outside the discriminant locus $\Sigma_s$ of $\phi\circ p_s$. It is
the set of $x\in X$ where $s_{|Y_x }\in \Gamma (Y_x , L^k)$ is not
transverse to the zero section. We set $X^0:=X-\Sigma_s$ and
$\tly_s^0:=(\phi\circ p_s)^{-1}(X-\Sigma_s)$ so that $\phi\circ p~: {\tly_s^0}
\lra X^0$ becomes a smooth family.  Then, since $\nab^0$ vanishes, the
above theorem implies the following

\begin{cor} \label{Gr2}
The vector bundle $\eb^0_{prim}=\eb^0=\rc^{r} (\phi\circ
p_s)_\star\oc_{\tly_s^0/X^0}= \oplus_{i=0}^{k-1}\rc ^{r} \phi_\star
\left( L^{-i}\right)$ with the Hodge metric is
Griffiths (actually even Nakano) semi-negative.
\end{cor}

\subsection{Singularities of the Hodge metric}

We now deal with the general case, namely $f : Y \lra B$ may not be
smooth.  We first recall the method of Fujita~\cite{fu}. We therefore
assume that the base $B$ is one dimensional.

The Hodge metric on the direct image of the relative canonical sheaf
is described as follows.  Let $b \in B$ be a point and let $(U, t)$ be
a local coordinate centered at $b = \{t= 0\}$.  A section $\om \in
\Gamma (U, f_\star K_{Y/B})$ -- when regarded as a section in $\Gamma
(U, f_\star Hom (f^\star K_B, K_Y))$ and applied to $f^\star dt$ --
gives a section of $K_Y$ on $f^{-1}(U)$ which we denote by $\om \cdot
dt$.  If $\varphi_b\in\Gamma (Y_b,K_{Y_b})$ fulfills the relation $\om
\cdot dt=\varphi_b\wedge f^\star dt$ over $Y_b$ (which amounts to
saying that in the differentiable trivialization $Y_{|U}\simeq
Y_b\times U$, the section $\om$ is sent to $\varphi_b$), then the
Hodge norm at $b \in B$ of the section $\om$ is
$$ \n \om \n^2_{Hodge} = (\sqrt{-1})^n (-1)^{\frac{n(n-1)}{2}}
        \int_{Y_b}\varphi_b\wedge\overline{\varphi_b},
$$ here $n = \dim Y -1$.  Fujita checked that in this setting in case
$\dim B =1$ the Hodge metric on $f_\star K_{Y/B}$ is bounded from
below by a positive quantity and hence that the only possible
singularities of the Hodge metric on $f_\star K_{Y/B}$ are poles
(see~\cite[lemma 1.12]{fu}).

We just give the typical example which occurs for a local model of our
cyclic covers (for some positive integer $m$).
$$
\begin{array}{ccccc}
\tly_s = \{(t,z,l)\in\cb^3 / l^k=t+z^m\} & \stackrel{p_s}\lra &
        Y&\stackrel{\phi}\lra&X = \{ t \in \cb \}\\ (t,z,l) &
        \mapsto&(t,z) & \mapsto & t
\end{array}
$$
The cotangent bundle $\Omega^1_{Y_s}$ is generated by $dt, dz, dl$
subject to the relation $kl^{k-1}dl-dt-mz^{m-1}dz=0$.  If $\om \cdot
dt$ is written as $\eta(z,l) dz\wedge dl$ for a holomorphic function
$\eta(z,l)$, then $\varphi_0$ may be chosen to be $\varphi_0 = k^{-1}
l^{1-k} \eta(z,l) dz$ with a pole of order $k-1$ on the fiber over
$t=0$ and no singularities elsewhere.

We can now go back to the proof of our main lemma.
\begin{proof}[Proof of lemma~\ref{mainlemma}]
Let $f~:~Z^{n+r}\lra X^n$ be a proper K{\"a}hler surjective map between
two complex manifolds of relative dimension $r$.  Let
$z=(z_1,z_2,\cdots, z_{n+r})$ and $x=(x_1,x_2,\cdots,x_n)$ be local
holomorphic coordinates on $Z$ and $X$ around $z_0$ and $x_0:=f(z_0)$.
The map $f$ is locally given by $(f_i(z))_{1\leq i\leq n}$. Let
$\Psi\in f_\star (K_{Z/X})$ and write
\begin{eqnarray*}
\Psi\cdot f^\star (dx_1\wedge dx_2\wedge\cdots\wedge dx_n) &=&\psi
dz_1\wedge dz_2\wedge\cdots\wedge dz_{n+r}\\ &=&\psi \det\left(
\frac{\d f_i}{\d z_j}\right)_{1\leq i\leq n\atop 1\leq j\leq n}
^{-1}f^\star (dx_1\wedge dx_2\wedge\cdots\wedge dx_n) \wedge
dz_{n+1}\wedge\cdots\wedge dz_{n+r}
\end{eqnarray*}
the Hodge norm of $\Psi$ is hence
\begin{eqnarray}\label{hodge.metric}
\n \Psi\n_{Hodge}^2=\int_{Z_x} |\psi|^2 \left|\det\left( \frac{\d
f_i}{\d z_j}\right) _{1\leq i\leq n\atop 1\leq j\leq n}\right|^{-2}
(\sqrt{-1})^r dz_{n+1}\wedge d\overline{z_{n+1}}\cdots\wedge dz_{n+r}
\wedge d\overline{z_{n+r}}.
\end{eqnarray}
The Hodge metric can hence only acquire poles, located furthermore
over the discriminant locus of $f$.  If now $\Psi$ belongs to $f_\star
(Jac_f\otimes K_{Z/X})$ then $\psi$ is a combination $\displaystyle
\psi=\sum %_{ J\subset \{ 1,2,\cdots , n+r\} \atop length(J)=n}
\psi_J\det\left( \frac{\d f_i}{\d z_j}\right)_{1\leq i\leq n\atop j\in
  J}$,
where the sum is taken over all multi-indexes $J \subset \{ 1,2,\cdots
, n+r\}$ of length $|J| =n$, manely $J = \{j_1, j_2, \ldots, j_n\}$
with $ 1 \leq j_1 < j_2 < \ldots < j_n \leq n+r$.  It follows that
$$ \Psi \cdot f^\star (dx_1\wedge dx_2\wedge\cdots\wedge dx_n)
=f^\star (dx_1\wedge dx_2\wedge\cdots\wedge dx_n)\wedge \bigg(\sum
%_{ J\subset \{ 1,2,\cdots , n+r\}, length(J)=n
%\atop M=\{1,2,\cdots , n+r\}_J} 
\psi_J dz_{m_1}\wedge \cdots\wedge dz_{m_r}\bigg).
$$ Here the last sum is taken over all multi-indexes $J \subset \{
1,2,\cdots , n+r\}$ of $|J| =n$, and $\{m_1, \ldots, m_r\} = \{
1,2,\cdots , n+r\} \setminus J$.
\end{proof}

In our setting, we take a global section $s$ of $L^k$ transverse to
the zero section and consider the corresponding cyclic cover $\tly_s$
of $Y$. Note that if the section $s$ is locally given by $\sigma
(x_1,x_2,\cdots, x_n,y_1,y_2,\cdots,y_r)$ and the map $\phi$ by
$(x_1,x_2,\cdots, x_n, y_1,y_2,\cdots,y_r)\mapsto (x_1,x_2,\cdots,
x_n)$, then there is an $j$ such that $\displaystyle \frac{\d s}{\d
y_j}(p_0)\not =0$ and the functions $(x_1,x_2,\cdots,x_n,
y_1,y_2,\cdots, y_{j-1},t,y_{j+1},\cdots,y_r)$ can be chosen as
coordinates on $\tly_s$ near $y_0\in p_s^{-1}(y_0)$, or there is an
$i$ such that $\displaystyle\frac{\d s}{\d x_i}(p_0)\not =0$ and the
coordinates on $\tly_s$ near $y_0$ can be chosen to be
$(x_1,x_2,\cdots x_{i-1},t,x_{i+1},\cdots ,x_n,y_1,y_2, \cdots ,
y_r)$.  In the first case, $\phi\circ p_s$ is a submersion at $y_0$
and its Jacobian ideal is $\oc_{\tly_s}$ at $y_0$. In the second case
the Jacobian ideal contains
\begin{eqnarray}\label{sing} 
\frac{\d x_i}{\d t}=\frac{\d s}{\d t} \left( \frac{\d s}{\d
x_i}\right)^{-1}\sim t^{k-1}.
\end{eqnarray} In
any case it contains $\displaystyle p_s^\star \oc
(-\frac{k-1}{k}D_s)$.  This has to be compared with the formula
$K_{Y_s/Y}=p_s^\star ((k-1)L)$.
%%%%%%%%%%%%%%%%%%%%%%%%%%%%%%%%%%%%%%%%%%%%%%%%%%%%%%
\section{Proof of Theorems}

\begin{proof}[Proof of Theorem~\ref{semi-ample}]
Let $\phi~:~Y\lra X$ be a projective submersion of relative dimension
$r$ between two complex manifolds and $L$ be a semi-ample line bundle
on $Y$.  Note that the relative projectivity condition is only used to
ensure the local freeness of $\phi_\star(K_{Y/X}\otimes L)$.

We take a positive integer $k$ so that $L^k$ is generated by its
global sections. Take a global section $s$ of $L^k$ transverse to the
zero section and consider the corresponding cyclic cover $\tly_s$ of
$Y$.
% and rank $S^k E > \dim X$.
We have seen that the Hodge metric on $\rc^r (\phi\circ p_s)_\star
\oc_{\tly_s}$ is smooth (as a Hodge metric on a smooth family)
non-degenerate and semi-negatively curved outside the discriminant
locus $\Sigma_{s_\alpha}$, is continuous on $X$, and may acquire zeros
at the points $x$ over which the section $s_{|\pb (E_x)}$ is
identically zero.  We now explain how the semi-ampleness assumption on
the line bundle $L$ helps to remove those singularities of the Hodge
metric.

When $L$ is semi-ample, $D\in |L^k|$ moves.  We choose $n+r+1$ generic
sections $s_\alpha$ that generate $L^k$.  The metric $h=\sum_{\alpha =
1}^{n+r+1} h_\alpha$ on $\rc^r \phi_\star ( L^{-1})\subset \rc^r
(\phi\circ p_s)_\star \oc_{\tly_s}$ is continuous on $X$, smooth
outside $\Sigma_h:=\cup_\alpha \Sigma_{s_\alpha}$ and without zeros.

Next let us discuss its curvature property.  Recall the formula
\begin{eqnarray}
\label{formule}
\frac{\langle\Theta (E,h)\xi,\xi\rangle}{\n \xi\n ^2} &=&
-\sqrt{-1}\d\!\db\log\n \xi\n ^2 +\frac{\langle \nabla_{E,h}
\xi,\nabla_{E,h} \xi\rangle}{\n \xi\n ^2} -\frac{\sqrt{-1}\d \n \xi\n
^2\wedge\db\n \xi\n ^2}{\n \xi\n ^4}\\ \nonumber &\geq&
-\sqrt{-1}\d\!\db\log\n \xi\n ^2
\end{eqnarray}
for a nowhere zero local holomorphic section $\xi$ of a holomorphic
vector bundle $E$ equipped with a smooth hermitian metric $h$. The
last two terms give the norm at $x$ of the fundamental form of the
inclusion $\oc_X \xi\subset E$.  Take a point $x_0$ in $X-\Sigma_h$, a
non-zero vector $\xi_0 \in \rc^r \phi_\star ( L^{-1})_{x_0}^\star$,
and a nowhere zero local holomorphic section $\xi\in\Gamma (U,\rc^r
\phi_\star ( L^{-1}))$ achieving the value $\xi_0$ at $x_0$ and normal
at $x_0$ for the metric $h$ (i.e. $\nabla_h \xi(x_0)=0$).  Then, the
last two terms in the formula~(\ref{formule}) vanish at $x_0$.
The corollary~\ref{Gr2} now asserts the function
$\log\n\xi\n_{h_\alpha} ^2$ -- whose complex Hessian (or Levi form) is
the opposite of the curvature of a line sub-bundle of $\rc^r
\phi_\star ( L^{-1})$ -- is plurisubharmonic on $U$.  It then follows
that $\log\n \xi\n^2_h = \log(\sum_\alpha \n\xi\n_{h_\alpha} ^2)$ is
plurisubharmonic also on $U$. This gives the Griffiths semi-negativity
of $\rc^r \phi_\star ( L^{-1})$ on $X -\Sigma_h$.  Since $h$ is
continuous and $\Sigma_h$ is an analytic subset (of zero Lebesgue
measure), we can conclude that the continuous metric $h$ on $\rc^r
\phi_\star ( L^{-1})$ is Griffiths semi-negative on the whole of $X$.
\end{proof}

\begin{proof}[Proof of Theorem~\ref{ample}]
We now assume that $L$ is ample. Here, the ampleness assumption
ensures the vanishing of $\mathcal{R}^{1}\phi_\star(K_{Y/X}\otimes L)$
which in turn easily ensures the local freeness of
$\phi_\star(K_{Y/X}\otimes L)$.
 
Let us recall Legendre-type formula applied for a metric $h =
\sum_{\alpha = 1}^\ell h_\alpha $ on $\rc^r \phi_\star ( L^{-1})$
gotten from different cyclic coverings $Y_{s_\alpha}$ (for a
$(1,0)$-form $u$, $|u|^2$ denotes $\sqrt{-1}u\wedge\overline{u}$)~:
$$ \sqrt{-1}\d\!\db \log(\sum_\alpha \n \xi\n_{h_\alpha} ^2) =
\frac{\sum_\alpha \n \xi\n_{h_\alpha} ^2 \sqrt{-1}\d\!\db\log\n
\xi\n_{h_\alpha} ^2} {\sum_\alpha \n \xi\n_{h_\alpha} ^2} +
\frac{\sum_{\alpha<\beta} \left|\d\log\n \xi\n_{h_\alpha} ^2 -
\d\log\n \xi\n_{h_\beta} ^2\right|^2 \n \xi\n_{h_\alpha}^2\n
\xi\n_{h_\beta}^2} {(\sum_\alpha \n \xi\n_{h_\alpha} ^2)^2}.
$$ Applying Griffiths curvature formula (Corollary~\ref{Gr2}) for
individual covering and the formula (\ref{formule}) for a line
sub-bundle we infer that in the right hand side, the first term is
semi-positive.  We need to add further Hodge metrics $h_\alpha$ to
make the second term -- hence the left hand side -- strictly positive.

The explicit expression in formula~(\ref{hodge.metric})
and~(\ref{sing}) will help to translate the algebraic ampleness
assumption on $L$ into a negativity property for a well chosen metric
on $\rc^r \phi_\star ( L^{-1})$.  Hence, we have for $\Psi \in
(\phi\circ p_s)_\star (K_{Y_s/Y})$
\begin{eqnarray*}
\d\log h_s(\Psi) &=& \d\log\int_{(\phi\circ p_s)^{-1}(x)} |\psi|^2
 \left|\left( \frac{\d s}{\d t} \right)^{-1} \frac{\d s}{\d
 x_i}\right|^{2} (\sqrt{-1})^r dz_{n+1}\wedge
 d\overline{z_{n+1}}\cdots\wedge dz_{n+r} \wedge d\overline{z_{n+r}}\\
 &=&\d\log\int_{Y_x } |\psi|^2 |s|^{-2k+2} \left|\frac{1}{k}\frac{\d
 s}{\d x_i}\right|^{2} (\sqrt{-1})^r dz_{n+1}\wedge
 d\overline{z_{n+1}}\cdots\wedge dz_{n+r} \wedge d\overline{z_{n+r}}
\end{eqnarray*}
We now take a positive integer~$k$ so large that the map
\begin{eqnarray*}
H^0(Y,L^k)&\lra&
H^0\left(Y,L^k\otimes\phi^\star(\oc_X/\mc^3_{x})\right)
\end{eqnarray*}
is surjective for every $x \in X$.  By the compactness of $Y$, we can
henceforth choose enough, but a finite number of sections $s_\alpha
\in H^0(Y,L^k)$ to ensure positivity in all the directions in the
Legendre formula. This gives continuous hermitian metrics on
$\phi_\star(K_{Y/X}\otimes L)$ with Griffiths positive
curvature. Using a regularization process as described in~\cite{mo}
these metrics may be smoothed keeping Griffiths positivity of the
curvature.
\end{proof}

\begin{rem}
Griffiths \cite[proposition 2.16]{gr} showed that the operator
$\nab^p~:\eb^p\lra\Omega^1_B\otimes\eb^{p-1}$ can be expressed as a
cup product with the Kodaira-Spencer class $\rho\in
\Omega^1_{B,b}\otimes H^1(Y_b,TY_b)$ of the family $f : Y \lra B$
coupled with a natural pairing.  In our setting, this in turn can be
related with the infinitesimal displacement of the hypersurfaces
$D_{s,x}$ of $Y_x $ given by the vanishing of the section $s_{|Y_x }$,
namely (see~\cite[chapter 5.2 (c)]{kodaira})
$$\begin{array}{cccccc} TX&\lra&H^0(D_{s,x}, \oc
(D_{s,x})_{|D_{s,x}})&\Big( \stackrel{\delta^\star}\lra &H^1(D_{s,x},
TD_{s,x})&\Big)\\ v&\mapsto& (\partial_{\widetilde{v}}
s)_{|D_{s,x}}&\Big(\mapsto&\rho (v)&\Big)
\end{array}$$
where $\widetilde{v}$ is a holomorphic vector field lifting $v$ on
$Y_x$. The map $\delta^\star$ is the co-boundary map in the long exact
sequence associated with the short exact sequence for the normal
bundle of the divisor $D_{s,x}$
$$0\lra TD_{s,x}\lra (TY_x)_{|D_{s,x}}\lra \oc
(D_{s,x})_{|D_{s,x}}\lra 0.$$ Our computations make explicit the idea
that $L$ being ample, the sections $s_{|Y_x}$ move sufficiently to
make the operator $\nab^p~:\eb^p\lra \Omega^1_B\otimes\eb^{p-1}$ have
non-zero contribution in the curvature formula.
\end{rem}

%%%%%%%%%%%%%%%%%%%%%%%%%%%%%%%%%%%%%%%%%%%%%%%%%%%%%%

\bigskip

Christophe Mourougane / Institut de Math{\'e}matiques de Jussieu / Plateau
7D / 175, rue du Chevaleret~/ 75013 Paris / France.\\ \textit{email} :
\textrm{mourouga@math.jussieu.fr}

\bigskip

Shigeharu Takayama / Graduate School of Mathematical Sciences /
University of Tokyo / 3-8-1 Komaba / 153-8914 / Japan.\\
\textit{email} : \textrm{taka@ms.u-tokyo.ac.jp}

%%%%%%%%%%%%%%%%%%%%%%%%%%%%%%%%%%%%%%%%%%%%%%%%%%%%%%

%%% Local Variables: 
%%% mode: latex
%%% TeX-master: t
%%% End: 
\end{document}